\documentclass{amsart}[12pt]
\usepackage{amssymb,amsmath}
\usepackage{enumerate}
\usepackage[english]{babel}
\usepackage{color}

\newcommand{\Om} {\Omega}
\newcommand {\ep} {\varepsilon}
\newcommand {\om} {\omega}

\newcommand {\ii} {\infty}

\newcommand {\lb} {\lambda}

\newcommand {\sm} {\setminus}
\newcommand {\su} {\subset}

\newcommand {\pr} {\prime}
\newcommand {\mc} {\mathcal}

\newcommand {\mbb} {\mathbb}

\newtheorem{teo}{Theorem}[section]

\theoremstyle{definition}

\newtheorem{df}{Definition}[section]

\DeclareMathOperator*{\esssup}{ess\,sup}

\title{On Egorov's Theorem for Infinite Measure}
\begin{document}
\date{August 18, 2023}

\begin{abstract}
A simple proof of Egorov's theorem for infinite measure is given.
\end{abstract}

\author{Semyon Litvinov}
\address{76 University Drive, Pennsylvania State University, Hazleton 18202}
\email{snl2@psu.edu}

\maketitle

\section{Introduction}
Let $\Om=(\Om,\Sigma,\mu)$ be a measure space, and let $\mc L^0(\Om)$ be the space of real-valued measurable functions on $\Om$.  Let $\|\cdot\|_\ii$ stand for the uniform norm in the space $\mc L^\ii(\Om)\su\mc L^0(\Om)$ of essentially bounded functions, that is, given $f\in\mc L^\ii(\Om)$,
\[
\|f\|_\ii=\esssup_{\om\in\Om}|f(\om)|.
\]
Denote by $\chi_E$ the characteristic function of a set $E\in\Sigma$.

\vskip5pt
\begin{df}
Let $\{f_n\}_{n=1}^\ii\su\mc L^0(\Om)$.

(A.E.) The sequence $\{f_n\}$ is said to converge to a function $f\in\mc L^0(\Om)$ {\it almost everywhere} (a.e.) if there exists a set $E\in\Sigma$ such that 
\[
\mu(\Om\sm E)=0\text{ \  and \ } f_n(\om)\to f(\om)\ \ \forall\ \om\in E.
\]

(A.U.) The sequence $\{f_n\}$ is said to converge to a function $f\in\mc L^0(\Om)$ {\it almost uniformly} (a.u.) if for every $\ep>0$ there exists a set $E\in\Sigma$ such that 
\[
\mu(\Om\sm E)<\ep\text{ \ and \ }\|(f_n-f)\chi_E\|_\ii\to 0.
\]
\end{df}
\vskip2pt
The classical Egorov's theorem can be stated as follows - see, for example, \cite{rf}:
\begin{teo}\label{t1} Assume that $\mu(\Om)<\ii$.
If $\{f_n\}\su\mc L^0(\Om)$ is such that $f_n\to f$ a.e., then $f_n\to f$ a.u.
\end{teo}
\noindent
It is easy to see that if $\mu(\Om)=\ii$, Theorem \ref{t1} is no longer valid: consider the space $\Om=(0,\ii)$ equipped with Lebesgue measure and let $f_n=\chi_{(n-1,n]}$, or $f_n=\chi_{[n,\ii)}$.
\vskip5pt

In \cite{ba}, the author, among other results, found a necessary and sufficient condition - which is in essence condition (\ref{e1}) below - under which an a.e. convergent sequence in an infinite measure space converges a.u. Note that the argument in \cite{ba} does not employ Theorem \ref{t1} - it rather follows "the steps ordinarily employed in measure theory to prove Egorov's and related theorems." The goal of this note is to provide a simple argument - based on an application of Egorov's theorem - to prove its infinite-measure extension:

\begin{teo}\label{t2}
Assume that $\mu(\Om)=\ii$. If $\{f_n\}\su\mc L^0(\Om)$ is such that $f_n\to f$ a.e., then $f_n\to f$ a.u. if and only if 
\begin{equation}\label{e1}
\forall \ \lb>0\ \ \exists \ n_\lb\text{ \ such that \ } \mu\Big\{\sup_{n\ge n_\lb}|f_n-f|\ge\lb\Big\}<\ii.
\end{equation}
\end{teo}

\section{Proof of Theorem \ref{t2}}
Without loss of generality, assume that $f=0$ a.e. 
\vskip3pt
If there exists $\lb>0$ such that $\mu\big\{\sup_{n\ge m}|f_n|\ge\lb\big\}=\ii$ for all $m\in\mbb N$, then $f_n\to 0$ a.u. clearly fails, so the "only if" part follows.
\vskip3pt
Next, assuming that (\ref{e1}) holds with $f=0$ a.e., we need to show that $f_n\to 0$ a.u. Fix $\ep>0$ and, given $m\in\mbb N$, let $n_1(m)$ be such that $\mu(F_m)<\ii$, where
\[
F_m=\Big\{\sup_{n\ge n_m}|f_n|\ge\frac1m\Big\},\text{\ so that \ } \big\||f_n|\chi_{\Om\sm F_m}\big\|_\ii\leq\frac1m\ \ \forall \ n\ge n_1(m).
\]
Then Theorem \ref{t1} implies that there exists\, $\Sigma\ni E_m\su F_m$ such that
\[
\mu(F_m\sm E_m)\leq\frac\ep{2^m}\text{\, \ and\, \ } \|f_n\chi_{E_m}\|_\ii\to 0\text{ \ as\ } n\to\ii.
\]
As the sequence $\{f_n\chi_{E_m}\}$ is Cauchy in $\mc L^\ii(\Om)$ for each $m$, it follows that 
\[
\forall\, m\ \ \exists\ n_2(m)\in\mbb N\text{ \ such that \ } \|(f_{n^\pr}-f_{n^{\pr\pr}})\chi_{E_m}\|_\ii\leq\frac1m\ \ \forall \ n^\pr, n^{\pr\pr}\ge n_2(m).
\]

Now, setting $G_m=E_m\cup (\Om\sm F_m)$ and $n(m)=\max\{n_1(m),n_2(m)\}$, we have $\mu(\Om\sm G_m)\leq\displaystyle\frac\ep{2^m}$ and
\[
\begin{split}
\|(f_{n^\pr}-f_{n^{\pr\pr}})\chi_{G_m}\|_\ii&=\|(f_{n^\pr}-f_{n^{\pr\pr}})\chi_{E_m}\|_\ii+\|(f_{n^\pr}-f_{n^{\pr\pr}})\chi_{\Om\sm F_m}\|_\ii\\
&\leq\frac1m+\big\||f_{n^\pr}|\chi_{\Om\sm F_m}\big\|_\ii+\big\||f_{n^{\pr\pr}}|\chi_{\Om\sm F_m}\big\|_\ii\leq\frac3m
\end{split}
\]
for each $m$ and $n^\pr, n^{\pr\pr}\ge n(m)$.
\vskip5pt
\noindent
Furthermore, letting\, $G=\bigcap_mG_m$ yields
\[
\mu(\Om\sm G)\leq\ep\text{\,\ \ and\,\ \ } \|(f_{n^\pr}-f_{n^{\pr\pr}})\chi_G\|_\ii\leq\frac3m \ \ \forall \ n^\pr, n^{\pr\pr}\ge n(m),
\]
implying that the sequence $\{f_n\chi_G\}$ is also Cauchy in $\mc L^\ii(\Om)$. 
\vskip5pt
Therefore, there exists $f\in\mc L^\ii(\Om)$ such that $\|f_n\chi_G-f\|_\ii\to 0$, so we have 
\[
\forall\ \ep>0 \ \ \exists\ G\in\Sigma\text{ \ such that\,\ } \mu(\Om\sm G)\leq\ep\text{ \ and \ }\|(f_n-f)\chi_G\|_\ii\to 0.
\]
Thus $f_n\to f$ a.u., so $f_n\to f$ a.e., which, due to $f_n\to 0$\, a.e., entails $f=0$\, a.e.
\vskip5pt
Summarizing, we conclude that 
\[
\forall\ \ep>0 \ \ \exists\ G\in\Sigma\text{ \ such that\,\ } \mu(\Om\sm G)\leq\ep\text{ \ and \ }\|f_n\chi_G\|_\ii\to 0,
\]
that is, $f_n\to 0$ a.u., thus completing the argument.

\end{document}